\documentclass[journal]{IEEEtran}
\usepackage{graphicx}
\usepackage{float}
\usepackage{amsthm}
\usepackage{amssymb}
\usepackage{textcomp}
\usepackage{nomencl}
\makenomenclature
\usepackage{ifthen}
\usepackage{setspace}
\usepackage{multirow}
\usepackage{paralist}
\usepackage[table,xcdraw]{xcolor}

\usepackage{cite}
\usepackage[pagebackref=false,breaklinks=true,letterpaper=true,colorlinks,bookmarks=false,citecolor=blue,linkcolor=blue]{hyperref}
\usepackage{amsmath}
\usepackage{amsmath,bm}
\usepackage{algorithmic}
\usepackage{algorithm}

\usepackage{array}

\ifCLASSOPTIONcompsoc
  \usepackage[caption=false,font=normalsize,labelfont=sf,textfont=sf]{subfig}
\else
  \usepackage[caption=false,font=footnotesize]{subfig}
\fi
\begin{document}
%
\title{A Mixed-Integer SDP Solution Approach to Distributionally Robust Unit Commitment with Second Order Moment Constraints}



\author{Xiaodong~Zheng,~\IEEEmembership{Student~Member,~IEEE,}
        Haoyong~Chen,~\IEEEmembership{Senior~Member,~IEEE,}
        Yan~Xu,~\IEEEmembership{Senior~Member,~IEEE,}
        Zhengmao~Li,~\IEEEmembership{Student~Member,~IEEE,}
        Zhenjia Lin,
        Zipeng Liang~\IEEEmembership{Student~Member,~IEEE}
\thanks{X. Zheng, H. Chen, Z. Lin and Z. Liang are with the School of Electric Power, South China University of Technology, Guangzhou 510641, China. X. Zheng is also with the School of Electrical and Electronic Engineering, Nanyang Technological University, Singapore (e-mail: z.xiaodong@mail.scut.edu.cn; eehychen@scut.edu.cn; epjack.lin@mail.scut.edu.cn; liangzipeng.ye@163.com).}
\thanks{Y. Xu and Z. Li are with the School of Electrical and Electronic Engineering, Nanyang Technological University, Singapore (email: eeyanxu@gmail.com; lizh0049@e.ntu.edu.sg).}}

\markboth{ }
{Shell \MakeLowercase{\textit{et al.}}: Bare Demo of IEEEtran.cls for IEEE Journals}
%



\IEEEtitleabstractindextext{
\begin{abstract}
A power system unit commitment (UC) problem considering uncertainties of renewable energy sources is investigated in this paper, through a distributionally robust optimization approach. We assume that the first and second order moments of stochastic parameters can be inferred from historical data, and then employed to model the set of probability distributions. The resulting problem is a two-stage distributionally robust unit commitment with second order moment constraints, and we show that it can be recast as a mixed-integer semidefinite programming (MI-SDP) with finite constraints. The solution algorithm of the problem comprises solving a series of relaxed MI-SDPs and a subroutine of feasibility checking and vertex generation. Based on the verification of strong duality of the semidefinite programming (SDP) problems, we propose a cutting plane algorithm for solving the MI-SDPs; we also introduce a SDP relaxation for the feasibility checking problem, which is an intractable biconvex optimization. Experimental results on a IEEE 6-bus system are presented, showing that without any tunings of parameters, the real-time operation cost of distributionally robust UC method outperforms those of deterministic UC and two-stage robust UC methods in general, and our method also enjoys higher reliability of dispatch operation.
\end{abstract}

\begin{IEEEkeywords}
Distributionally robust optimization, mixed-integer semidefinite programming, probability distribution, renewable energy integration, unit commitment.
\end{IEEEkeywords}}

\maketitle

\IEEEdisplaynontitleabstractindextext

%
\IEEEpeerreviewmaketitle



\section{Introduction}
%
%
%
%
\IEEEPARstart{W}{ith} the characteristics of sustainability, environmental friendliness and economic efficiency, renewable energy sources (RES) are being increasingly integrated into the power systems worldwide \cite{mathiesen2011100, dai2016green}. One common concern of developing power systems with high share of RES is how to cope with their intermittence and fluctuation natures, i.e., how to ensure the reliability and efficiency of system operation under the uncertainties of renewable energy generation \cite{purvins2011challenges, badal2019survey}.

In the recent decades, an enormous literature has been devoted to improving the reliability and efficiency of short-term operation, by using cutting-edge optimization techniques such as stochastic programming and robust optimization, etc. Among these methods, stochastic programming method relies on an accurate probability distribution of stochastic parameters \cite{shapiro2000stochastic}, since otherwise the decisions yielded may be misleading. Another method, robust optimization attracts more attentions in recent years especially on the topic of unit commitment (UC) \cite{jiang2012robust, zeng2013solving, bertsimas2013adaptive, zheng2019hierarchical}. Instead of minimizing the expectation of operation cost, robust UC minimizes the worst-case total cost over all possible realizations within a pre-defined uncertainty set, and hence has less demand on the information of renewable energy generation.

For the robust UC method, the decision quality depends on some components of the model, such as the uncertainty set. If the uncertainty set is too conservative, it includes some extreme scenarios that are unlikely to realize. {\color{black}These extreme scenarios (as well as other scenarios within the uncertainty set), however, are supposed to be attainable with probability one. As a result, the solution may commit redundant units and enforce too much reserve.} Actually, in the era of big-data \cite{chen2018learning}, it is reasonable to make use of available data to improve the quality of decision-making. For example, we can draw an empirical probability distribution of wind power (or its prediction error) from historical data. As the data increase in size, the empirical distribution should be close enough to the ``true'' distribution. Therefore, we can expect a better solution through optimizing the (worst-case) expected total cost under a \emph{set of distributions} centered at the empirical one, which leads to the concept of distributionally robust optimization.

Given the popularity and strong modeling capacity of distributionally robust optimization, we will give a detailed review of this method, starting from general problems, and then shifting to those applied into the UC problem.

One of the pioneer work of distributionally robust optimization is \cite{delage2010distributionally}, in which the authors deal with the single-stage distributionally robust optimization problem under moment uncertainty. The mean (first order moment) of stochastic parameters is forced to lie in an ellipsoid centered at the estimate, while the covariance matrix (\emph{second order moment}) is forced to lie in the intersection of two positive semi-definite cones. The model is recast into a SDP with infinite number of constraints indexed by the stochastic parameters, and it is further shown to be solvable (under some mild assumptions) in time polynomial in the dimensions of decision variables and stochastic parameters. As for the two-stage setup (which is more practical in power system operations), the authors in \cite{bertsimas2010models} propose a SDP model for a class of two-stage distributionally robust optimization (``minimax two-stage stochastic linear optimization problems with risk aversion'' by the authors) with known first and second order moments. The authors prove that for the problem with random objective, a tight SDP formulation can be achieved; for the problem with random right-hand side, it is $\mathcal{NP}$-hard in general.
Beyond moment constraints methods, another effective modeling technique is the distance-based method, which enforces the distributions within a fixed distance from the empirical one based on a certain metric. For example, \cite{zhao2018data} studies the two-stage distributionally robust optimization (``data-driven risk-averse stochastic optimization'' by the authors) with Wasserstein metric, and shows that it can be successfully reformulated to the traditional two-stage robust optimization problem. Notice that these works haven't considered discrete decision variables.

In recent years, the problem of unit commitment with distributional uncertainty has been addressed by increasing number of papers, such as \cite{liu2015stochastic, xiong2017distributionally, zhao2018distributionally, duan2018data, ding2019duality, chen2018distributionally}. As one of the pioneers of distributionally robust unit commitment, \cite{xiong2017distributionally} uses linear moment constraints and linear decision rule (LDR) to gain a tractable formulation. However, there may be a natural defect caused by LDR as will be discussed later. In \cite{zhao2018distributionally}, the authors deal with the probability of outage for devices (generation units and transmission lines), and the set of distributions is defined by the expected number of failure devices, which is also a linear first order moment constraint. \cite{duan2018data} models the set of distributions based on the confidence bands for cumulative distribution function (CDF), and incorporates dtributionally robust chance constraints into the second-stage dispatch problem as well. In \cite{ding2019duality}, the set of wind power distribution is restricted by the 1-norm and $\infty$-norm with respect to the empirical probability, and a pre-designed confidence level is necessary to design this set. The inner max-min problem can be decomposed into independent dispatch problems to be solved in parallel. In \cite{chen2018distributionally}, another type of metric, i.e., Kullback-Leibler divergence, is employed to model the set of distributions for the UC problem.

However, the abovementioned works seldom incorporate second order moment, and thus the variance and correlation information of/between stochastic parameters is discarded. In \cite{wei2016distributionally}, a two-stage distributionally robust optimization model with first and second order moments is developed to deal with the uncertainties in a joint energy and reserve dispatch problem. Since discrete variables don't exist in the dispatch problem, the optimization model is reformulated as a SDP with finite constrains indexed by the number of vertexes of a polyhedral. A two-step solution algorithm is proposed and the effectiveness of the distributionally robust optimization model is verified on the IEEE 118-bus system.

The principle goal of this paper is to deploy the second order moment constraints into the distributionally robust UC problem, and develop an efficient solution algorithm. The contributions of this paper include:
\begin{enumerate}
  \item The distributionally robust UC model with second order moment constraints is developed and then successfully recast as a MI-SDP problem.
  \item A cutting plane algorithm is proposed to solve the MI-SDP problems, with optimality and convergence guarantee provided by a relevant proposition.
  \item In the feasibility checking subroutine of the solution algorithm, a SDP relaxed problem is formulated to derive the lower bound of the biconvex feasibility checking problem, and thus guide the subroutine.
  \item A case study based on a one-year wind power dataset, comprising 30 days' optimal scheduling and simulation, is carried out. Thorough numerical results and comparisons with those of deterministic UC and robust UC methods are delivered.
\end{enumerate}

The remainder of the paper is organized as follows: Section II describes the formulation of the second-order moment constrained distributionally robust UC problem. Section III introduces the reformulation procedure and a detailed solution algorithm. Section IV presents the case study and numerical results. Section V concludes with discussion.


\section{Problem Formulation}
The full UC model is removed to Appendix~\ref{append:UC}, and more detailed models can be found in literature, e.g., \cite{chen2016key, bertsimas2013adaptive}. Hereafter, we work on the compact formulation below,
\begin{subequations} \label{eqn:UC}
    \begin{align}
     \nonumber \underset{\bm{x,y}}{\mathop{\min }}\,&\bm{c}_{I}^{\top}\bm{x}+\bm{c}_{C}^{\top}\bm{y}\\
     \operatorname{s.t.~}&\bm{Fx}\ge \bm{f} \label{eqn:UC_Ff} \\
     & \bm{Gy}\ge \bm{g} \label{eqn:UC_Gg} \\
     & \bm{Ax}+\bm{By}\ge \bm{d} \label{eqn:UC_ABd} \\
     & \bm{Uy}=\bm{Ww}, \label{eqn:UC_Uw}
    \end{align}
\end{subequations}
where $\bm{x}$, $\bm{y}$ are binary unit commitment variables and continuous dispatch variables (including variables of RES spillage and load shedding) respectively, and $\bm{c}_I$, $\bm{c}_C$ are the cost coefficients associated with them. Constraint~(\ref{eqn:UC_Ff}) denotes the state transition equations of units and minimum on/off time limits of units; constraint~(\ref{eqn:UC_Gg}) includes the ramping limits of units and power flow limits of transmission lines; the coupling constraint~(\ref{eqn:UC_ABd}) represents the generation capacity constraints, etc.; constraint~(\ref{eqn:UC_Uw}) is the power balance condition, e.g., the DC power flow equation.

To account for the uncertainties of nodal injections brought by intermittent renewable energy sources, a uncertainty/stochastic parameter $\bm{\xi} \in \Xi $ is integrated into (\ref{eqn:UC_Uw}), yielding $\bm{Uy}+\bm{V\xi} =\bm{Ww}$. Now we can write out the robust counterpart of (\ref{eqn:UC}) as a two-stage optimization problem,
\begin{align}
\nonumber & \underset{\bm{x}}{\mathop{\min }}\,\underset{\bm{\xi} \in \Xi  }{\mathop{\max }}\,\bm{c}_{I}^{\top}\bm{x}+Q(\bm{x,\xi} ) \\
& \operatorname{s.t.~}\bm{Fx}\ge \bm{f}
\label{eqn:RUC}
\end{align}
where
\begin{align}
\nonumber Q(\bm{x},&\bm{\xi} )=\underset{\bm{y}}{\mathop{\min }}\, \bm{c}_{C}^{\top}\bm{y} \\
\nonumber \operatorname{s.t.~}& \bm{Gy}\ge \bm{g} \\
\nonumber & \bm{Ax}+\bm{By}\ge \bm{d} \\
& \bm{Uy}+\bm{V\xi} =\bm{Ww}.
\label{eqn:Q}
\end{align}
It is shown in equation~(\ref{eqn:RUC}) that robust UC seeks the worst-case scenario within an uncertainty  set and derives an optimal UC solution under that case. The most commonly used uncertainty set in power system optimizations is the following polyhedral set \cite{chen2016key, bertsimas2013adaptive, jiang2012robust, zheng2019hierarchical},
\begin{equation}
\begin{aligned}
 \Xi=\left\{ \bm{\xi} \left| {{\xi}_{i}}\in  \right.\left[ {{\bar{\xi}}_{i}}-k\sigma_i,{{\bar{\xi}}_{i}}+k\sigma_i \right],~
\sum\nolimits_{i}{\left| {{\xi}_{i}}-{{\bar{\xi}}_{i}} \right|/k\sigma_i\le {\Gamma}} \right\}
\end{aligned}
\label{eqn:uncertainty_set}	
\end{equation}
where $\sigma_i$ is the standard deviation (std.) of the $i$-th parameter of $\bm{\xi}$, $k$ is determined by the chosen confidence level, and $\Gamma$ is the uncertainty budget.

In distributionally robust optimizations, however, one would like to find a worst-case probability distribution within a set of distributions. Building upon the spirit of distributionally robust optimization, we present a distributionally robust UC model (DRUC) as follows,
\begin{align}
\nonumber & \underset{\bm{x}}{\mathop{\min }}\,\underset{\mathbb{P}\in \mathsf{\mathcal{P}}}{\mathop{\max }}\,\bm{c}_{I}^{\top}\bm{x}+{{\mathbb{E}}_{\mathbb{P}}}[Q(\bm{x,\xi} )] \\
& \operatorname{s.t.~} \bm{Fx}\ge \bm{f}
\label{eqn:DRUC}
\end{align}
where $Q(\bm{x,\xi})$ is the value function of $\bm{x}$ and $\bm{\xi}$ as defined in equation~(\ref{eqn:Q}), and $\mathsf{\mathcal{P}}$ is the set of distributions defined as,
\begin{equation} \label{eqn:set_P}
\mathsf{\mathcal{P}}=\left\{ \mathbb{P}\in {{\mathsf{\mathcal{P}}}_{0}}(\Xi  )\left|
\begin{aligned}
  & {{\mathbb{E}}_{\mathbb{P}}}[1]=1 \\
 & {{\mathbb{E}}_{\mathbb{P}}}[\bm{\xi} ]=\bar{\bm{\xi} } \\
 & {{\mathbb{E}}_{\mathbb{P}}}[\bm{\xi} {{\bm{\xi} }^{\top}}]=\bm{\Sigma} +\bar{\bm{\xi} }{{\bar{\bm{\xi} }}^{\top}} \\
\end{aligned}
\right. \right\}.
\end{equation}

In the definition of $\mathsf{\mathcal{P}}$, ${{\mathsf{\mathcal{P}}}_{0}}(\Xi  )$ denotes the set of all of probability measures on a sigma algebra of $\Xi  \subseteq {{\mathbb{R}}^{\left| \bm{\xi}  \right|}}$, and $\bm{\Sigma} \in {{\mathbb{R}}^{\left| \bm{\xi}  \right|\times \left| \bm{\xi}  \right|}}$ is the covariance matrix. {\color{black}The last constraint in Eq.~(\ref{eqn:set_P}) indicates that the covariance of random parameters is $\bm{\Sigma}$, i.e., $\mathbb{E}[(\bm{\xi}-\bar{\bm{\xi}})(\bm{\xi}-\bar{\bm{\xi}})^{\top}] = \bm{\Sigma}$.} Note that set $\mathsf{\mathcal{P}}$ is shaped by the exact first order moment $\bar{\bm{\xi}}$ and second order moment $\bm{\Sigma}$, which are the a-priori estimates of stochastic parameters obtained from statistics analysis of historical data. Given a set of historical data $\left\{ {{\bm{\xi} }_{i}}\left| i=1,...,M \right. \right\}$, many methods can be applied to estimate the parameters of its probability distribution, including moment methods, maximum likelihood estimation methods, and Bayesian nonparametric models \cite{wang2016wind, campbell2015bayesian}, etc. In this paper, we adopt a simple unbiased moment estimator \cite{delage2010distributionally}:
\begin{equation}
  \bar{\bm{\xi} }=\frac{1}{M}\sum\limits_{i=1}^{M}{{{\bm{\xi} }_{i}}},~~\bm{\Sigma} =\frac{1}{M-1}\sum\limits_{k=1}^{M}{\left( {{\bm{\xi} }_{i}}-\bar{\bm{\xi} } \right)}{{\left( {{\bm{\xi} }_{i}}-\bar{\bm{\xi} } \right)}^{\top }}.
\end{equation}

Although $\bm{\Sigma}$ is a positive semidefinite matrix by construction\footnote{{\color{black}For any vector $\bm{z}\in\mathbb{R}^{|\bm{\xi}|}$, we have $\bm{z}^{\top}\bm{\Sigma}\bm{z} = \bm{z}^{\top}\mathbb{E}[(\bm{\xi}-\bar{\bm{\xi}})(\bm{\xi}-\bar{\bm{\xi}})^{\top}]\bm{z} = \mathbb{E}[\bm{z}^{\top}(\bm{\xi}-\bar{\bm{\xi}})(\bm{\xi}-\bar{\bm{\xi}})^{\top}\bm{z}] =  \mathbb{E}[(\bm{z}^{\top}(\bm{\xi}-\bar{\bm{\xi}}))^2] \ge 0 $. Therefore, a covariance matrix $\bm{\Sigma}$ must be positive semidefinite, i.e., $\bm{\Sigma} \succeq 0$.}}, we assume that $\bm{\Sigma} \succ 0$ (in the case we study, we check this condition and it always holds). The necessity of this assumption will be clear shortly.

\section{Solution Method}

\subsection{Reformulation as a Mixed-Integer SDP}
In the distributionally robust UC model (\ref{eqn:DRUC}), the second-stage maximizing problem is an infinite dimensional conic linear problem, whose explicit formulation is,
\begin{align}\label{eqn:Z_primal}
\nonumber  Z(\bm{x})& =\underset{\mathbb{P}\in \mathsf{\mathcal{P}}}{\mathop{\max }}\,{{\mathbb{E}}_{\mathrm{}}}[Q(\bm{x,\xi} )] \\
\nonumber & =\underset{{{f}_{\bm{\xi} }}}{\mathop{\max }}\,\int_{\Xi  }{Q(\bm{x,\xi} ){{f}_{\bm{\xi} }}\mathrm{d}\bm{\xi} } \\
\nonumber \operatorname{s.t.~}& {{f}_{\bm{\xi} }}\ge 0 ~~\forall \bm{\xi} \in \Xi   \\
\nonumber & \int_{\Xi }{{{f}_{\bm{\xi}}}\mathrm{d}\bm{\xi}=1}~~:{{h}_{0}} \\
\nonumber & \int_{\Xi }{{{\xi}_{i}}{{f}_{\bm{\xi}}}\mathrm{d}\bm{\xi}={{\bar{ \xi }}_{i}}} ~~ :{{h}_{i}},\text{ }i=1,2,...,\left| \bm{\xi} \right| \\
\nonumber & \int_{\Xi }{{{\xi}_{j}}{{ \xi }_{k}}{{f}_{\bm{\xi}}}\mathrm{d}\bm{\xi}={{\Sigma}_{jk}}+{{\bar{ \xi }}_{j}}{{\bar{\xi}}_{k}}} ~~ :{{H}_{jk}}, \\
& ~~~~~~~~~~~~~j,k=1,2,...,\left| \bm{\xi} \right|,j\ge k.
\end{align}
Problem (\ref{eqn:Z_primal}) has infinite number of decision variables ${{f}_{\bm{\xi} }}$ indexed by $\bm{\xi}$, and finite number of moment constraints, so its dual problem may be more tractable. In the spirit of linear programming duality theory, we associate a dual variable with each equality constraint{\color{black}, and dualizing problem~(\ref{eqn:Z_primal}) into a minimization problem. The process of dualization is detailed in Appendix~\ref{append:dualization}}. According to Theorem 16.3.2 in \cite{bertsimas2000moment}, if the right-hand side moment vector lies in the interior of the set of feasible moments, then strong duality holds. It then concludes that $\bm{\Sigma} \succ 0$ is a sufficient condition for the holding of strong duality. Therefore, $Z(\bm{x})$ can be {\color{black}equivalently} represented by its dual problem,
\begin{align}\label{eqn:Z_dual}
\nonumber   & Z(\bm{x})=\underset{\bm{H,h},{{h}_{0}}}{\mathop{\min }}\,\operatorname{tr}\left( {{\bm{H}}^{\top}}\bm{\Theta}  \right)+{{\bm{h}}^{\top}}\bar{\bm{\xi} }+{{h}_{0}} \\
 & \operatorname{s.t.~}{{\bm{\xi} }^{\top}}\bm{H\xi} +{{\bm{h}}^{\top}}\bm{\xi} +{{h}_{0}}\ge Q(\bm{x,\xi} )~~ \forall \bm{\xi} \in \Xi,
\end{align}
where $\bm{H}$, $\bm{h}$, ${{h}_{0}}$ are dual variables associated with the equality constraints of problem~(\ref{eqn:Z_primal}), $\bm{\Theta} =\bm{\Sigma} +\bar{\bm{\xi} }{{\bar{\bm{\xi} }}^{\top}}$, and $\operatorname{tr}(\bm{X})$ denotes the trace of matrix $\bm{X}$ which returns the summation of the elements on the main diagonal.

It seems that problem (\ref{eqn:Z_dual}) or (\ref{eqn:DRUC}) can be recast as a single-stage robust optimization by introducing LDR to $Q(\bm{x,\xi} )$ \cite{xiong2017distributionally}. Yet, we found that LDR makes problem (\ref{eqn:Z_dual}) \emph{trivial}. To prove this, we can look back to the objective function of its primal problem (\ref{eqn:Z_primal}): with LDR, we have $Q(\bm{x,\xi} )=\bm{c}_{C}^{\top}{{\bm{y}}^{\text{opt}}}=\bm{c}_{C}^{\top}({{\bm{y}}^{0}}+\bm{R\xi} )$, where ${{\bm{y}}^{0}}$ is the optimal dispatch decision under the predicted scenario and $\bm{R}$ is the optimal LDR matrix. Since ${{\mathbb{E}}_{\mathbb{P}}}[\bm{\xi} ]=\bar{\bm{\xi} }$, we have,
{\color{black}
\begin{align}
\nonumber &\underset{\mathbb{P}\in \mathsf{\mathcal{P}}}{\mathop{\max }}\,{{\mathbb{E}}_{\mathbb{P}}}[Q(\bm{x,\xi} )]\\
\nonumber =&\underset{\mathbb{P}\in \mathsf{\mathcal{P}}}{\mathop{\max }}\,{{\mathbb{E}}_{\mathbb{P}}}[\bm{c}_{C}^{\top}({{\bm{y}}^{0}}+\bm{R\xi} )]\\
=&\bm{c}_{C}^{\top}({{\bm{y}}^{0}}+\bm{R}\bar{\bm{\xi} }),
\end{align}}
which shows that the optimum of problem~(\ref{eqn:Z_primal}) relies merely on $\bar{\bm{\xi}}$ but not on other information declared by set $\mathsf{\mathcal{P}}$, and hence the dual problem (\ref{eqn:Z_dual}) is also trivial.

By now, we realize that problem (\ref{eqn:Z_dual}) must be solved with the exact value of $Q(\bm{x,\xi} )$. As we know, $Q(\bm{x,\xi} )$ is the value function of $\bm{x}$ and $\bm{\xi}$, which can be represented by its dual problem, i.e.,
\begin{align}\label{eqn:Q_dual}
  \nonumber Q(\bm{x,\xi})=~& \underset{\bm{\lambda ,\mu ,\nu} }{\mathop{\max }}\,{{\bm{g}}^{\top}}\bm{\lambda} -{{(\bm{Ax}-\bm{d})}^{\top}}\bm{\mu} -{{(\bm{V\xi} -\bm{Ww})}^{\top}}\bm{\nu}  \\
 \nonumber \operatorname{s.t.~}& {{\bm{c}}_{C}}-({{\bm{G}}^{\top}}\bm{\lambda} +{{\bm{B}}^{\top}}\bm{\mu} +{{\bm{U}}^{\top}}\bm{\nu} )=\bm{0} \\
 & \bm{\lambda} \ge \bm{0},\bm{\mu} \ge \bm{0}.
\end{align}
Herein, $\bm{\lambda ,\mu ,\nu}$ are dual variables of problem (\ref{eqn:Q}). Problem (\ref{eqn:Q}) or (\ref{eqn:Q_dual}) is always feasible and bounded, because slack variables have been introduced into the equality constraint of (\ref{eqn:Q}). Then, as suggested by \cite{wei2016distributionally}, problem (\ref{eqn:Z_dual}) can be formulated as a standard optimization problem (\ref{eqn:Z_vert}) by listing all the vertexes of the feasible region $\Omega$ of (\ref{eqn:Q_dual}), i.e., $\text{vert}(\Omega  )=\left\{ ({{\bm{\lambda} }_{i}},{{\bm{\mu} }_{i}},{{\bm{\nu} }_{i}}) \right\}$,
\begin{subequations} \label{eqn:Z_vert}
    \begin{align}
    & Z(\bm{x})=\underset{\bm{H,h},{{h}_{0}}}{\mathop{\min }}\,\operatorname{tr}\left( {{\bm{H}}^{\top}}\bm{\Theta}  \right)+{{\bm{h}}^{\top}}\bar{\bm{\xi} }+{{h}_{0}} \\
    \nonumber \operatorname{s.t.~}& {{\bm{\xi} }^{\top}}\bm{H\xi} +{{\bm{h}}^{\top}}\bm{\xi} +{{h}_{0}}\ge {{\bm{g}}^{\top}}{{\bm{\lambda} }_{i}}-{{(\bm{Ax}-\bm{d})}^{\top}}{{\bm{\mu} }_{i}} \\
    & ~~~~~~~~~-{{(\bm{V\xi} -\bm{Ww})}^{\top}}{{\bm{\nu} }_{i}}~~\forall \bm{\xi} \in \Xi~~\forall i\in \mathsf{\mathcal{I}},   \label{eqn:Z_vert_cs}
    \end{align}
\end{subequations}
where $\mathsf{\mathcal{I}}$ is the set of indexes for $\text{vert}(\Omega  )$. Note that constraint (\ref{eqn:Z_vert_cs}) is equivalent with,
\begin{gather}\label{eqn:LMI}
    \begin{bmatrix} \bm{\xi}  \\ 1 \end{bmatrix}^{\top}
    \begin{bmatrix}
    \bm{H} &
    0.5(\bm{h}+{{\bm{V}}^{\top}}{{\bm{\nu} }_{i}}) \\
    0.5{{(\bm{h}+{{\bm{V}}^{\top}}{{\bm{\nu} }_{i}})}^{\top}} &
    \tilde{h}_{0}
    \end{bmatrix}
    \begin{bmatrix} \bm{\xi}  \\ 1 \end{bmatrix}
    \ge 0~\forall \bm{\xi} \in \Xi,
\end{gather}
where $\tilde{h}_{0}={{h}_{0}}-{{\bm{g}}^{\top}}{{\bm{\lambda} }_{i}}+{{(\bm{Ax}-\bm{d})}^{\top}}{{\bm{\mu} }_{i}}-{({\bm{Ww}})^{\top}}{{\bm{\nu} }_{i}}$. When $\Xi  ={{\mathbb{\bm{R}}}^{\left| \bm{\xi}  \right|}}$, condition~(\ref{eqn:LMI}) can be further represented by the following semidefinite cone,
\begin{gather}\label{eqn:SDP_cone}
    \begin{bmatrix}
    \bm{H} &
    0.5(\bm{h}+{{\bm{V}}^{\top}}{{\bm{\nu} }_{i}}) \\
    0.5{{(\bm{h}+{{\bm{V}}^{\top}}{{\bm{\nu} }_{i}})}^{\top}} &
    \tilde{h}_{0}
    \end{bmatrix}
    \succeq 0.
\end{gather}
Therefore, problem (\ref{eqn:DRUC}) equals to the following MI-SDP,
\begin{align}
\nonumber \underset{\bm{x,H,h},{{h}_{0}}}{\mathop{\min }}\,& \bm{c}_{I}^{\top}\bm{x}+\operatorname{tr}\left( {{\bm{H}}^{\top}}\bm{\Theta}  \right)+{{\bm{h}}^{\top}}\bar{\bm{\xi} }+{{h}_{0}} \\
\nonumber \operatorname{s.t.~}& \bm{Fx}\ge \bm{f} \\
 & (\ref{eqn:SDP_cone})~~\forall i\in \mathsf{\mathcal{I}}. \label{eqn:MISDP}
\end{align}

\subsection{Cutting Plane Algorithm}
MI-SDP has been investigated by literature such as \cite{gally2018framework, ni2018mixed}, and different solution methods have been proposed, including SDP-based Branch-and-Bound methods, outer approximation approaches, etc. However, since conic optimization cannot be solved very efficiently as other traditional convex problems nowadays, it is reasonable to develop a specialized algorithm for a certain MI-SDP problem in order to gain higher solution efficiency. In this paper, we propose a cutting plane method to solve problem (\ref{eqn:DRUC}).

Before going into the details of the cutting plane algorithm, we have the following Proposition:
\newtheorem{prop}{Proposition}
\begin{prop}
 Given any feasible moments $\bar{\bm{\xi} }$ and $\bm{\Sigma} \succ 0$, if for $\forall \bm{\xi} \in \Xi $, $Q(\bm{x,\xi} )$ exists and is bounded, then
   \begin{enumerate}
     \item problem (\ref{eqn:Z_primal}) is feasible and bounded; \label{prop:1_i}
     \item problem (\ref{eqn:Z_dual}) is feasible and bounded; \label{prop:1_ii}
     \item replace $\mathsf{\mathcal{I}}$ in MI-SDP problem (\ref{eqn:MISDP}) with any subset of $\mathsf{\mathcal{I}}$, denoted as ${{\mathsf{\mathcal{I}}}^{-}}$, the relaxed MI-SDP problem is also feasible and bounded; \label{prop:1_iii}
     \item replace $\mathsf{\mathcal{I}}$ with any subset ${{\mathsf{\mathcal{I}}}^{-}}$ and fix the first-stage variable $\bm{x}$ in problem (\ref{eqn:MISDP}), the resulting SDP given by (\ref{eqn:Z_vert_sub}) as well its dual problem (also a SDP) are both feasible, and moreover, strong duality holds. \label{prop:1_iv}
   \end{enumerate}
\label{prop:1}
\end{prop}

In Proposition \ref{prop:1}, \ref{prop:1_i}) and \ref{prop:1_ii}) are the bases of \ref{prop:1_iii}) and \ref{prop:1_iv}), whereas \ref{prop:1_iv}) {\color{black}certificates the strong duality and thus} leads to the cutting plane algorithm. Because at each iteration, the dual problem of (\ref{eqn:Z_vert_sub}) will be solved to generate a Benders cut, the cutting plane algorithm will converge to the \emph{global optimality} with probability one only if the problem (\ref{eqn:Z_vert_sub}) is convex and strong duality holds. {\color{black}It is noted that the convergence of Algorithm~\ref{alg:1} is basically guaranteed by the finiteness of the solution space of integer variables.}

Assuming at current step that a subset of $\text{vert}(\Omega )$ indexed by $\mathsf{\mathcal{I}}{{\mathsf{}}^{-}}$ is available, fix the binary variables and we will have the following SDP in hand,
\begin{align}\label{eqn:Z_vert_sub}
  \nonumber {{Z}^{-}}(\bm{x})=& \underset{\bm{H,h},{{h}_{0}}}{\mathop{\min }}\,\operatorname{tr}\left( {{\bm{H}}^{\top}}\bm{\Theta}  \right)+{{\bm{h}}^{\top}}\bar{\bm{\xi} }+{{h}_{0}} \\
  \operatorname{s.t.~}& (\ref{eqn:SDP_cone})~~ \forall i\in \mathsf{\mathcal{I^-}}.
\end{align}

It can be proved that the dual problem of (\ref{eqn:Z_vert_sub}) has the following formulation \cite{freund2004introduction},
\begin{align}\label{eqn:Z_vert_sub_dual}
  \nonumber &Z_{\mathrm{d}}^{-}(\bm{x})=\underset{\bm{\gamma ,\upsilon} }{\mathop{\max }}\,\sum\limits_{i\in \mathsf{\mathcal{I}}{{\mathsf{}}^{-}}}{\bm{b}_{i}^{\top}{{\bm{\gamma} }_{i}}} \\
 \nonumber &\operatorname{s.t.~} \tilde{\bm{C}}-\sum\limits_{j=1}^{\left| \bm{\xi}  \right|+1}{{{\gamma }_{i,j}}{{\tilde{\bm{A}}}_{j}}}+\sum\limits_{k=1}^{({{\left| \bm{\xi}  \right|}^{2}}+\left| \bm{\xi}  \right|)/2}{{{\upsilon }_{i,k}}{{\tilde{\bm{B}}}_{k}}}\succeq 0~~i=1 \\
 & -\sum\limits_{j=1}^{\left| \bm{\xi}  \right|+1}{{{\gamma }_{i,j}}{{\tilde{\bm{A}}}_{j}}}-\sum\limits_{k=1}^{({{\left| \bm{\xi}  \right|}^{2}}+\left| \bm{\xi}  \right|)/2}{{{\upsilon }_{i,k}}{{\tilde{\bm{B}}}_{k}}}\succeq 0~\forall i\in \mathsf{\mathcal{I}}{{\mathsf{}}^{-}}\backslash \left\{ 1 \right\}
\end{align}
where ${{\bm{\gamma} }_{i}}\in {{\mathbb{R}}^{\left| \bm{\xi}  \right|+1}}$ and ${{\bm{\upsilon} }_{i}}\in {{\mathbb{R}}^{({{\left| \bm{\xi}  \right|}^{2}}+\left| \bm{\xi}  \right|)/2}}$ are the dual variables for problem (\ref{eqn:Z_vert_sub}), ${{\gamma }_{i,j}}$ denotes the $j$-th element of the vector, the same goes for ${{\upsilon }_{i,k}}$, ${{\bm{b}}_{i}}={{\left( \begin{matrix}
   -0.5{{\bm{V}}^{\top}}{{\bm{\nu} }_{i}} & h_0-{{\widetilde{h}}_{0}}  \\
\end{matrix} \right)}}$ is the right-hand side of problem (\ref{eqn:Z_vert_sub}), and $\tilde{\bm{C}}$, ${\tilde{{\bm{A}}}_{j}}$, ${\tilde{{\bm{B}}}_{k}}$ are coefficient matrixes with proper dimensions and elements.

Now it is enough to present the cutting plane algorithm for solving (the relaxation of) MI-SDP (\ref{eqn:MISDP}):
\begin{algorithm}
\caption{Cutting Plane Algorithm}
\label{alg:1}
\begin{algorithmic}
\STATE \textbf{Initialize:} Let $L=1$, $\bm{x}_1$ be an initial UC solution, $UB= +\infty$, $LB= -\infty$. Choose a convergence tolerance $\epsilon$.
\WHILE {$UB-LB>\epsilon$}
\STATE $\bullet$ Solve the sub-problem $Z_{\mathrm{d}}^{-}({{\bm{x}}^{L}})$ to obtain $\bm{\gamma} _{i}^{L}$, set $UB=\bm{c}_{I}^{\top}{{\bm{x}}^{L}}+\sum_{i\in \mathsf{\mathcal{I}}{{\mathsf{}}^{-}}}{\bm{b}_{i}^{\top}\bm{\gamma} _{i}^{L}}$
\STATE $\bullet$ Solve the master problem
\begin{align*}
 \nonumber \underset{\bm{x},\theta}{\mathop{\min }}\, & \bm{c}_{I}^{\top}\bm{x}+\theta  \\
 \nonumber\operatorname{s.t.~} & \bm{Fx}\ge \bm{f} \\
 & \theta \ge \sum\limits_{i\in \mathsf{\mathcal{I}}{{\mathsf{}}^{-}}}{\bm{b}_{i}^{\top}\bm{\gamma} _{i}^{l}}~~l=1,...,L\text{,}
\end{align*}
record the optimal solution ${{\bm{x}}^{\text{*}}},{{\theta }^{*}}$, set ${{\bm{x}}^{L+1}}={{\bm{x}}^{\text{*}}}$, $LB=\bm{c}_{I}^{\top}{{\bm{x}}^{\text{*}}}+{{\theta }^{*}}$, $L\leftarrow L+1$.
\ENDWHILE
\STATE \textbf{return} $\bm{x}^L$
\end{algorithmic}
\end{algorithm}

\subsection{Feasibility Checking}
After solving the relaxed problem of (\ref{eqn:MISDP}), we should check whether the constraint of problem (\ref{eqn:Z_dual}) is violated for some new vertexes of $\Omega $, and add them to the relaxed MI-SDP problem if any. To make the narration more clear, let us rewrite the constraint of problem (\ref{eqn:Z_dual}):
\begin{align*}
{{\bm{\xi} }^{\top}}\bm{H\xi} +{{\bm{h}}^{\top}}\bm{\xi} +{{h}_{0}}\ge Q(\bm{x,\xi} )~~ \forall \bm{\xi} \in \Xi.
\end{align*}
As aforementioned, $Q(\bm{x,\xi} )$ can be represented by its dual problem, resulting in the constraint,
\begin{align}\label{eqn:Z_dual_cs}
  \nonumber &  {{\bm{\xi} }^{\top}}\bm{H\xi} +{{\bm{h}}^{\top}}\bm{\xi} +{{h}_{0}} \ge {{\bm{g}}^{\top}}\bm{\lambda} -{{(\bm{Ax}-\bm{d})}^{\top}}\bm{\mu} \\
 & -{{(\bm{V\xi} -\bm{Ww})}^{\top}}\bm{\nu}~~\forall \bm{\xi} \in \Xi ,~(\bm{\lambda ,\mu ,\nu} )\in \Omega.
\end{align}
After solving the relaxed problem of (\ref{eqn:MISDP}), we obtain the UC solution $\bm{x}$, and the solution values of $\bm{H, h}, h_0$. However, such solution are unnecessary within the feasible region depicted by constraint~(\ref{eqn:Z_dual_cs}). In \cite{wei2016distributionally}, the authors solve a biconvex problem (\ref{eqn:violation_check}) to check whether constraint (\ref{eqn:Z_dual_cs}) is violated for some $(\bm{\lambda ,\mu ,\nu})$, i.e., whether the objective value of (\ref{eqn:violation_check}) is negative. Since problem (\ref{eqn:violation_check}) is $\mathcal{NP}$-hard, the checking procedure is actualized through a sequential convex optimization approach as shown in \cite{wei2016distributionally}, which converges to a local optimum of (\ref{eqn:violation_check}) (hence the optimal solution provides an upper bound).
\begin{align}\label{eqn:violation_check}
  \nonumber \underset{\bm{\xi},\bm{\lambda ,\mu ,\nu}} {\mathop{\min}} & {{\bm{\xi} }^{\top}}\bm{H\xi} +{{\bm{h}}^{\top}}\bm{\xi} +{{h}_{0}}-{{\bm{g}}^{\top}}\bm{\lambda} +{{(\bm{Ax}-\bm{d})}^{\top}}\bm{\mu} \\
  \nonumber & +{{(\bm{V\xi} -\bm{Ww})}^{\top}}\bm{\nu}  \\
 \operatorname{s.t.~}& \bm{\xi} \in \Xi ,~(\bm{\lambda ,\mu ,\nu} )\in \Omega.
\end{align}

In what follows, we introduce a method to determine the lower bound of problem~(\ref{eqn:violation_check}). Note that the only nonconvex term of problem~(\ref{eqn:violation_check}) is ${{(\bm{V\xi} )}^{\top}}\bm{\nu}$, which can be written as,
\begin{gather}\label{eqn:QP_relaxation_cs}
{{ \begin{bmatrix}
   \bm{\xi}   \\
   \bm{\nu}   \\
\end{bmatrix} }^{\top}} \begin{bmatrix}
   0 & 0.5{{\bm{V}}^{\top}}  \\
   0.5{{\bm{V}}^{\top}} & 0  \\
\end{bmatrix}  \begin{bmatrix}
   \bm{\xi}   \\
   \bm{\nu}   \\
\end{bmatrix} =\operatorname{tr}(\widetilde{\bm{V}}\bm{\Psi} ),
\end{gather}
where
\begin{align*}
\widetilde{\bm{V}} = \left( \begin{matrix}
   0 & 0.5{{\bm{V}}^{\top}}  \\
   0.5{{\bm{V}}^{\top}} & 0  \\
\end{matrix} \right),
\bm{\Psi} =\bm{\psi} {{\bm{\psi} }^{\top}}=\left( \begin{matrix}
   \bm{\xi}   \\
   \bm{\nu}   \\
\end{matrix} \right){{\left( \begin{matrix}
   \bm{\xi}   \\
   \bm{\nu}   \\
\end{matrix} \right)}^{\top}}.
\end{align*}
Relax the constraint $\bm{\Psi} =\bm{\psi} {{\bm{\psi} }^{\top}}$ into $\bm{\Psi} \succeq \bm{\psi} {{\bm{\psi} }^{\top}}$ and then problem~(\ref{eqn:violation_check}) can be relaxed to a SDP,
\begin{align}\label{eqn:violation_check_ralaxation}
    \nonumber \underset{\bm{\xi},\bm{\lambda ,\mu ,\nu},\bm{\Psi,\psi}} {\mathop{\min}} & {{\bm{\xi} }^{\top}}\bm{H\xi} +{{\bm{h}}^{\top}}\bm{\xi} +{{h}_{0}}-{{\bm{g}}^{\top}}\bm{\lambda} +{{(\bm{Ax}-\bm{d})}^{\top}}\bm{\mu} \\
    \nonumber & -(\bm{Ww})^{\top}\bm{\nu} + \operatorname{tr}(\widetilde{\bm{V}}\bm{\Psi} )  \\
    \operatorname{s.t.~}& \bm{\xi} \in \Xi ,~(\bm{\lambda ,\mu ,\nu} )\in \Omega, ~\bm{\Psi} \succeq \bm{\psi} {{\bm{\psi} }^{\top}}.
\end{align}

The merit of problem (\ref{eqn:violation_check_ralaxation}) is that it provides a lower bound for problem (\ref{eqn:violation_check}), so that we can terminate the solution algorithm of the DRUC problem once the lower bound is nonnegative; or otherwise, repeat the sequential convex optimization approach from more initial points to explore new vertex, if the current upper bound is positive while the lower bound is a large negative number.

\section{Case Study}
\begin{figure}
  \centering
  \includegraphics[width=3in]{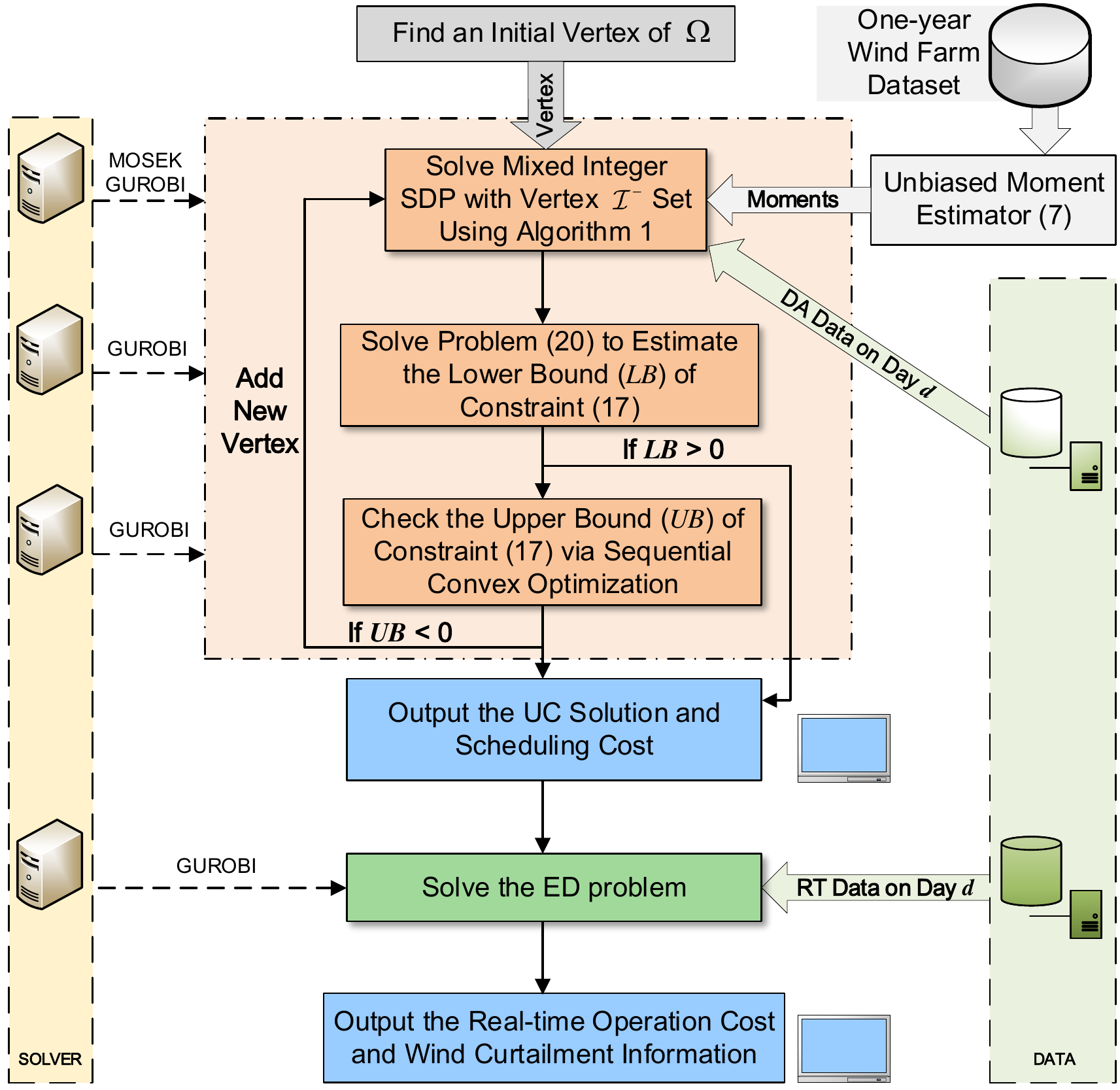}\\
  \caption{The framework of numerical experiments.}\label{fig:flowchart}
\end{figure}

In this section, numerical experiments are carried out to verify the effectiveness of the second order moment constrained DRUC model, and evaluate the performance of the solution approach. The whole framework of the experiments is illustrated by Fig. \ref{fig:flowchart}, which contains three main blocks, namely moments estimation of wind power prediction errors, solving the DRUC problem, and solving an economic dispatch (ED) problem to simulate the real-time operation. {\color{black}A penalty term has been augmented to the objective function to mimic the high cost of load shedding, wind curtailment or dispatching fast-start units; the penalty cost is set as \$1000/MWh.}

The optimization problems are built in YALMIP \cite{lofberg2004yalmip}, LPs, QPs and MILPs are solved by Gurobi 8.1.1\footnote{Available online: http://www.gurobi.com}, and SDPs are solved by MOSEK 8.1\footnote{Available online: http://www.mosek.com}. All runs are executed on an Intel i5 CPU machine running at 1.80 GHz with 8 GB of RAM.

\subsection{System Setup and Data}
\begin{figure}
  \centering
  \includegraphics[width=3in]{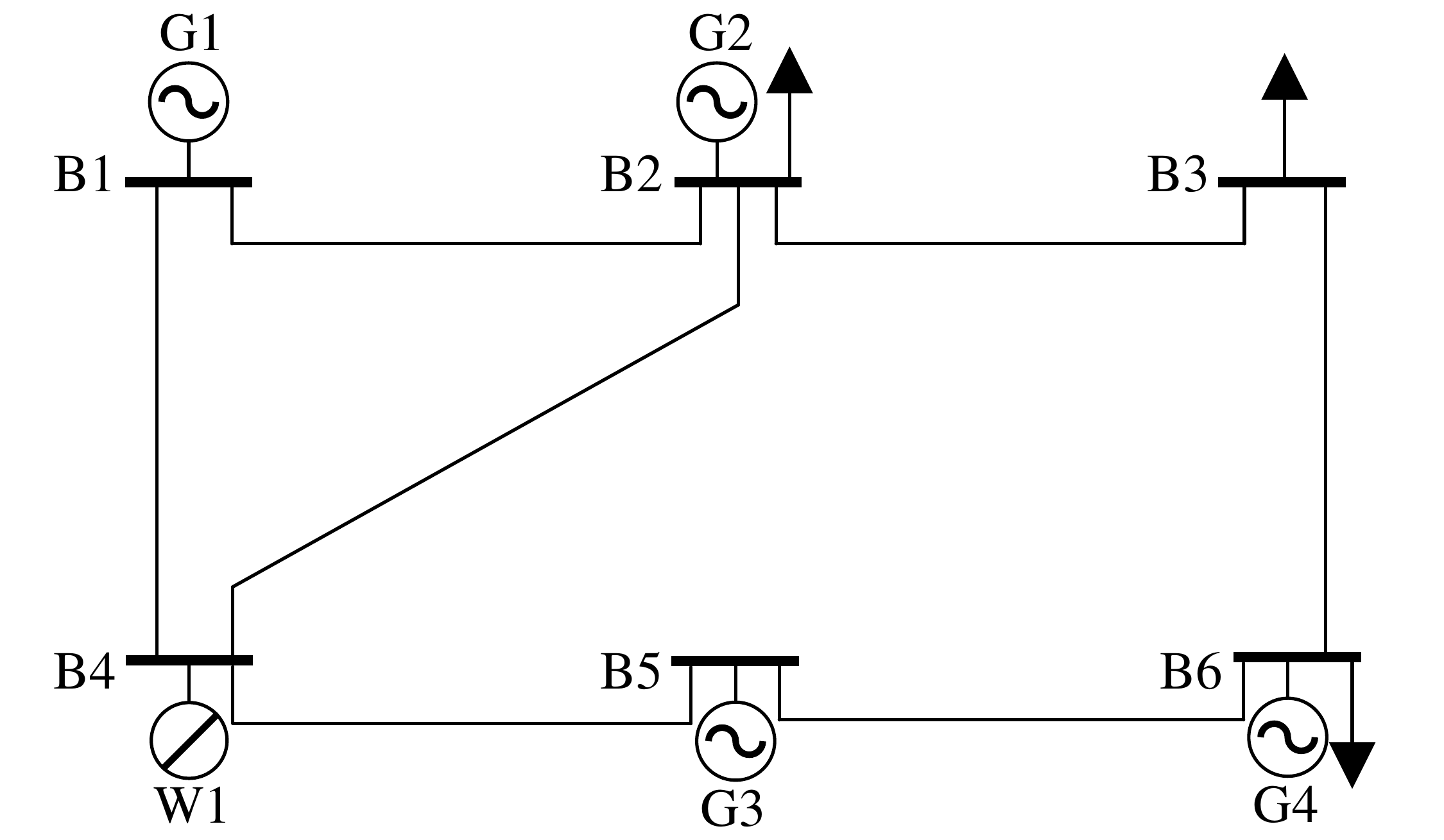}\\
  \caption{The standard IEEE 6-bus test system.}\label{fig:six_bus}
\end{figure}

The IEEE 6-bus system as shown in Fig. \ref{fig:six_bus} is used to do the experiments. The generator data, fuel/cost data, transmission line data as well as the locations of generator, load and wind farm are in line with those in \cite{jiang2012robust}. The system load demand information comes from a realistic power system, and the prediction error of load demand is neglected. The capacity of the wind farm installed at Bus 4 is 100 MW, namely 13.7\% of the total capacity of this system.

Time-synchronous hourly wind series\footnote{Both day-head predictions and real-time data are available. We choose the data of wind farm \#1, and scale them with an appropriate factor.} over one year provided by \cite{pena2018extended} are used for two purposes: \textit{i}) we feed the whole dataset over one year into the unbiased moment estimator to obtain the mean vector and covariance matrix of wind power prediction errors; \textit{ii}) we randomly draw out 30-days samples from the dataset, employ the day-ahead data to the scheduling models to derive unit commitment decisions, and the real-time data to run the simulation and evaluate the dispatch cost.

\begin{figure}[!t]
\centering
\subfloat[]{\includegraphics[width=3in]{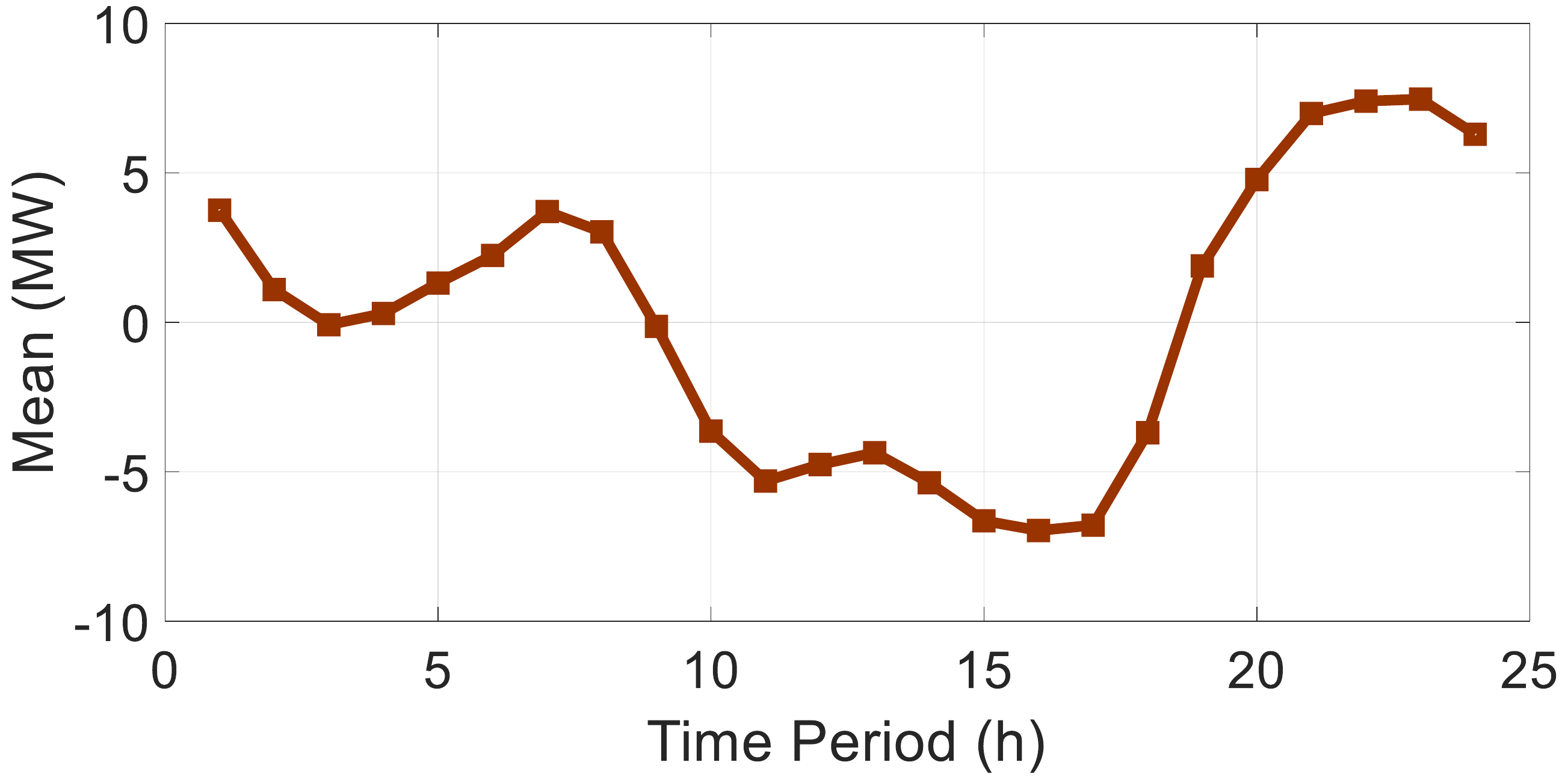}%
\label{fig:mean}}
\hfil
\subfloat[]{\includegraphics[width=3.3in]{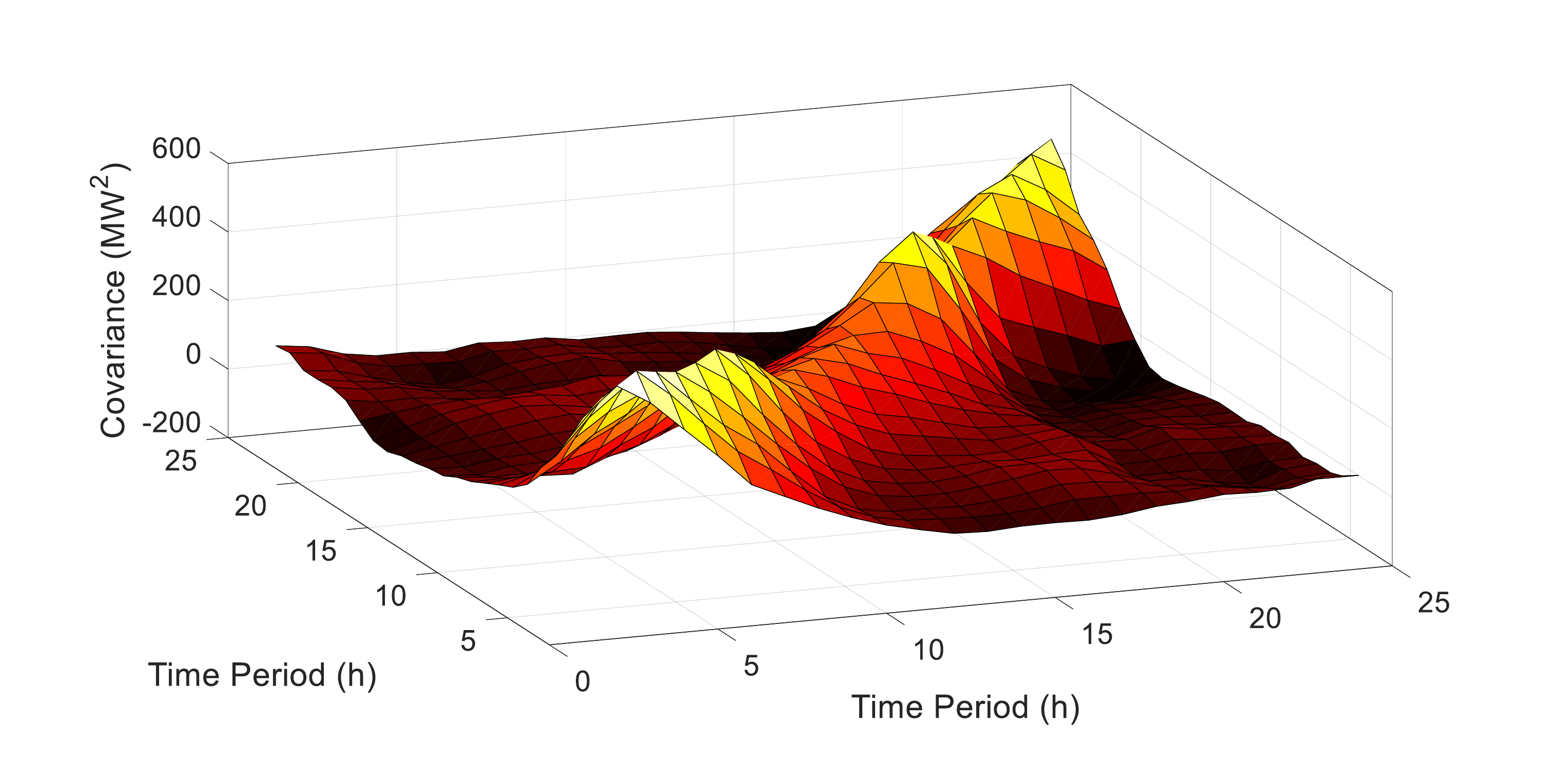}%
\label{fig:cov}}
\caption{The moment information of wind power prediction errors drawn from the one-year wind farm operation data. (a) Mean. (b) Covariance.}
\label{fig:mean_cov}
\end{figure}

The mean vector and covariance matrix of wind power prediction errors achieved by our statistics method are presented in Fig. \ref{fig:mean_cov}. It is shown in Fig. \ref{fig:mean_cov}\subref{fig:mean} that the means of prediction error are nonzero at almost all the scheduling period and also vary with hours. According to Fig. \ref{fig:mean_cov}\subref{fig:cov}, we know that the covariances of prediction errors between any two different time slots are relatively small compared with the variances (i.e., elements on the main diagonal), yet those between closer time slots are still considerable. The abovementioned information of wind power prediction errors can be explicitly modeled by set $\mathsf{\mathcal{P}}$, and hence not only the deviation rates but also the correlations of deviations will be accounted by the distributionally robust UC model.

\subsection{Cost Efficiency and Reliability of Dispatch Operation}
We compare our method (DRUC) with another two methods, namely deterministic UC (UC) and robust UC (RUC). For RUC, we use the conventional polyhedral uncertainty set (\ref{eqn:uncertainty_set}), in which the mean and the interval of uncertainty parameter are respectively set as $\bar{\bm{\xi}}$ and $\pm1.44\bm{\sigma}$ ($\bm{\sigma}$ is the vector of std. of $\bm{\xi}$). The solution approach for RUC contains the column-and-constraint generation method \cite{zeng2013solving} and the McCormick MILP reformualtion method \cite{chen2018novel}.

\begin{figure}[!t]
\centering
\subfloat[]{\includegraphics[width=2.8in]{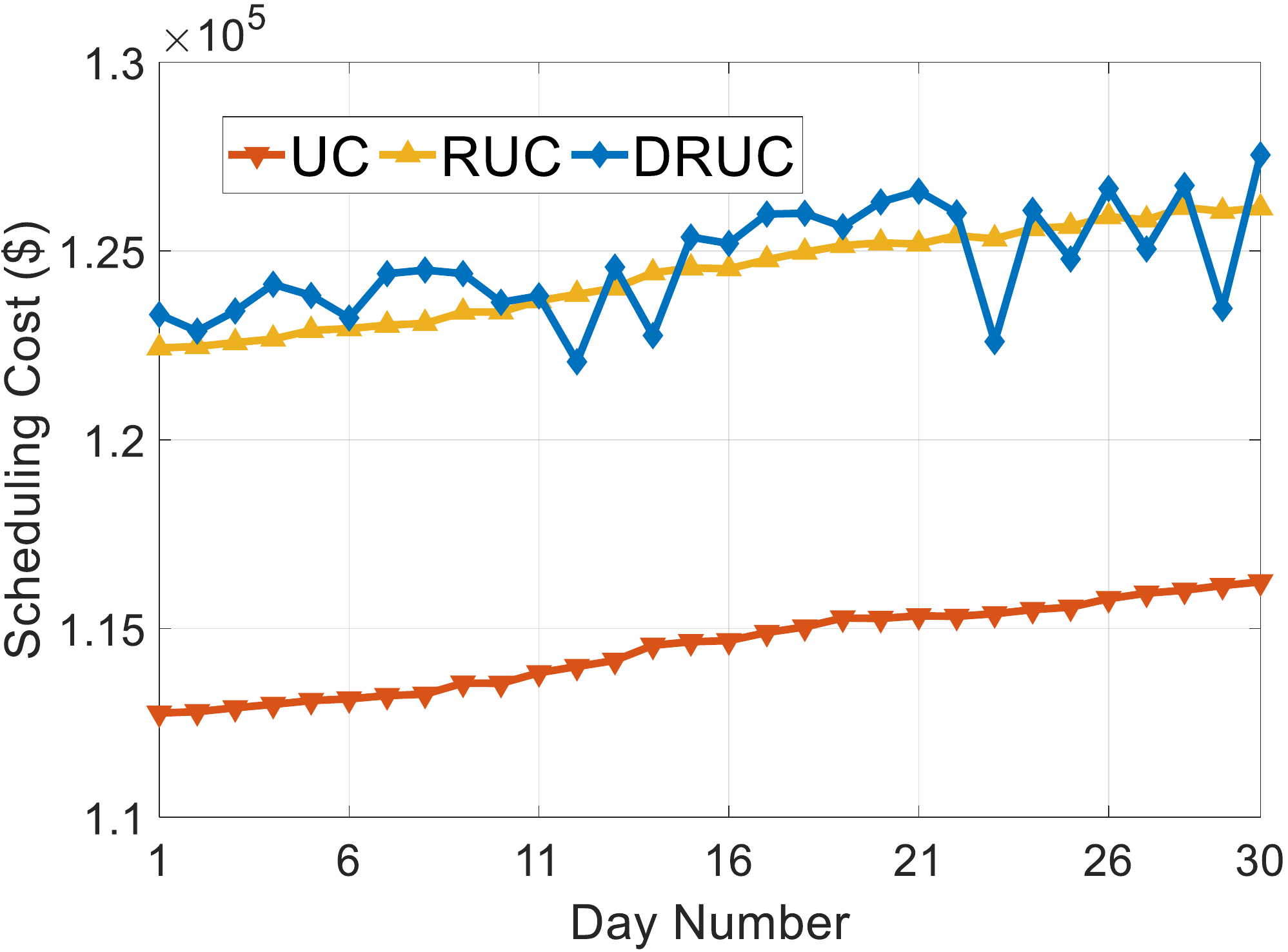}%
\label{fig:da_cost}}
\hfil
\subfloat[]{\includegraphics[width=2.8in]{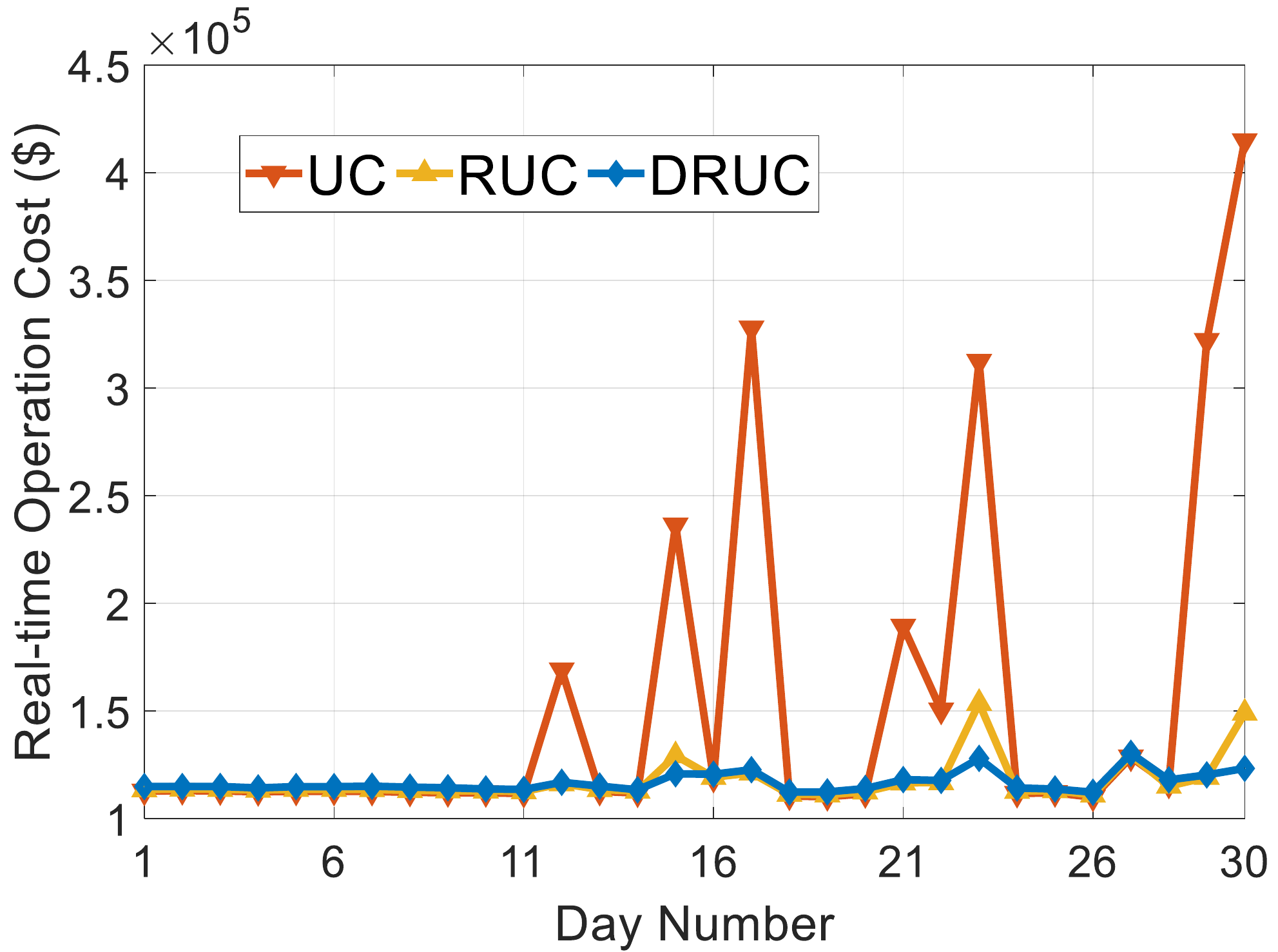}%
\label{fig:rt_cost}}
\caption{Comparisons of costs during 30 days. (a) Day-ahead scheduling costs. (b) Real-time operation costs.}
\label{fig:costs}
\end{figure}

To compare the cost efficiency, the day-ahead scheduling costs and real-time operation costs over 30 days are presented in Fig. \ref{fig:costs}. As shown is Fig. \ref{fig:costs}\subref{fig:da_cost}, the deterministic UC method always enjoys the lowest scheduling cost, because the model doesn't deal with uncertainties; both RUC and DRUC have higher scheduling costs, since the objective values of them are yielded from the worst-case scenario or probability distribution. As shown is Fig. \ref{fig:costs}\subref{fig:rt_cost}, in intra-day operation, the deterministic UC method suffers from wind curtailment on about 7 days, leading to extremely high costs; visible wind curtailment only appears on about 3 days for the RUC method; however, the {\color{black}real-time operation cost} of DRUC outperforms those of both UC and RUC, which is actually close to the scheduling cost. {\color{black}It is observed in Fig.~\ref{fig:costs}\subref{fig:rt_cost} that the real-time cost of the UC method can be several times higher than the scheduling cost. This is due to the high penalty cost and the size of the test system, i.e., if a smaller penalty cost is chosen, and the test system has a larger size and a higher demand level, then the total penalty cost will not dominate the generation cost.}

\begin{figure}
  \centering
  \small
  \includegraphics[width=2.8in]{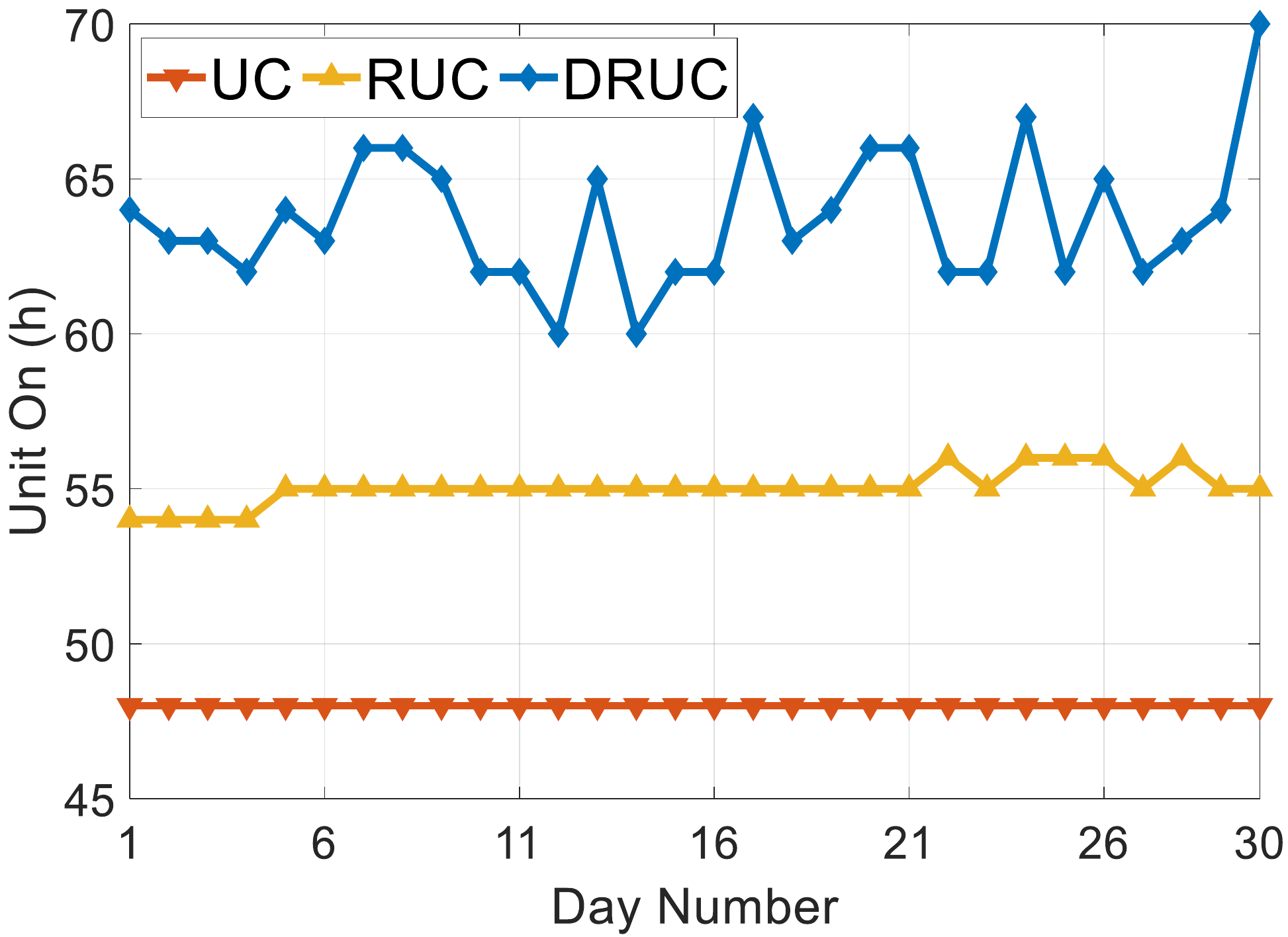}\\
  \caption{Number of on-hour for generation units during 30 days.}\label{fig:unit_on}
\end{figure}

In Table \ref{tab:results}, the statistics of mean and standard deviation of scheduling costs, {\color{black}real-time operation costs} and the amount of wind curtailment over 30 days are presented. We can see that the unit commitment derived from DRUC not only has highest cost efficiency (mean of {\color{black}real-time operation cost}), but also provides enough reliability (mean of wind spillage) and robustness (std. of both {\color{black}real-time operation cost} and wind spillage). In this case, the high efficiency and reliability of DRUC is achieved by committing more units, see Fig. \ref{fig:unit_on}. Such a reasonable decision is owed to the nature of distributionally robust optimization: in set $\mathsf{\mathcal{P}}$, the prediction errors of wind power are indicated by the first order moment constraint, while the correlation of prediction errors is further captured by the second order moment constraint. The more precise the uncertainty model is, the higher efficiency and reliability our decisions can achieve. On the contrary, it is quite challenging for robust models to account for the spatial and temporal correlations of uncertainties.
\begin{table}[t]
\centering{}%
\caption{Comparisons of Results between UC, RUC and DRUC methods.\label{tab:results}}
\footnotesize
\begin{tabular}{c|cccc}
\hline
\multicolumn{1}{l|}{}                   & \multicolumn{1}{l}{} & UC      & RSUC    & DRUC                           \\ \hline
                                        & mean                 & 114,495 & 124,375 & 124,696                        \\
\multirow{-2}{*}{Scheduling Cost (\$)} & std.                 & 1,146   & 1,249   & 1,435                          \\ \cline{2-5}
                                        & mean                 & 153,766 & 117,699 & {\color[HTML]{034e7b} 116,839} \\
\multirow{-2}{*}{Real-time Operation Cost (\$)}       & std.                 & 82,360  & 10,218  & 4,507                          \\ \cline{2-5}
                                        & mean                 & 39.94   & 2.40    & {\color[HTML]{034e7b} 0.33}    \\
\multirow{-2}{*}{Wind Spillage (MW)}    & std.                 & 81.00   & 7.00    & 1.00                              \\ \hline
\end{tabular}
\end{table}

\begin{figure}[!t]
\centering
\subfloat[]{\includegraphics[width=3.4in]{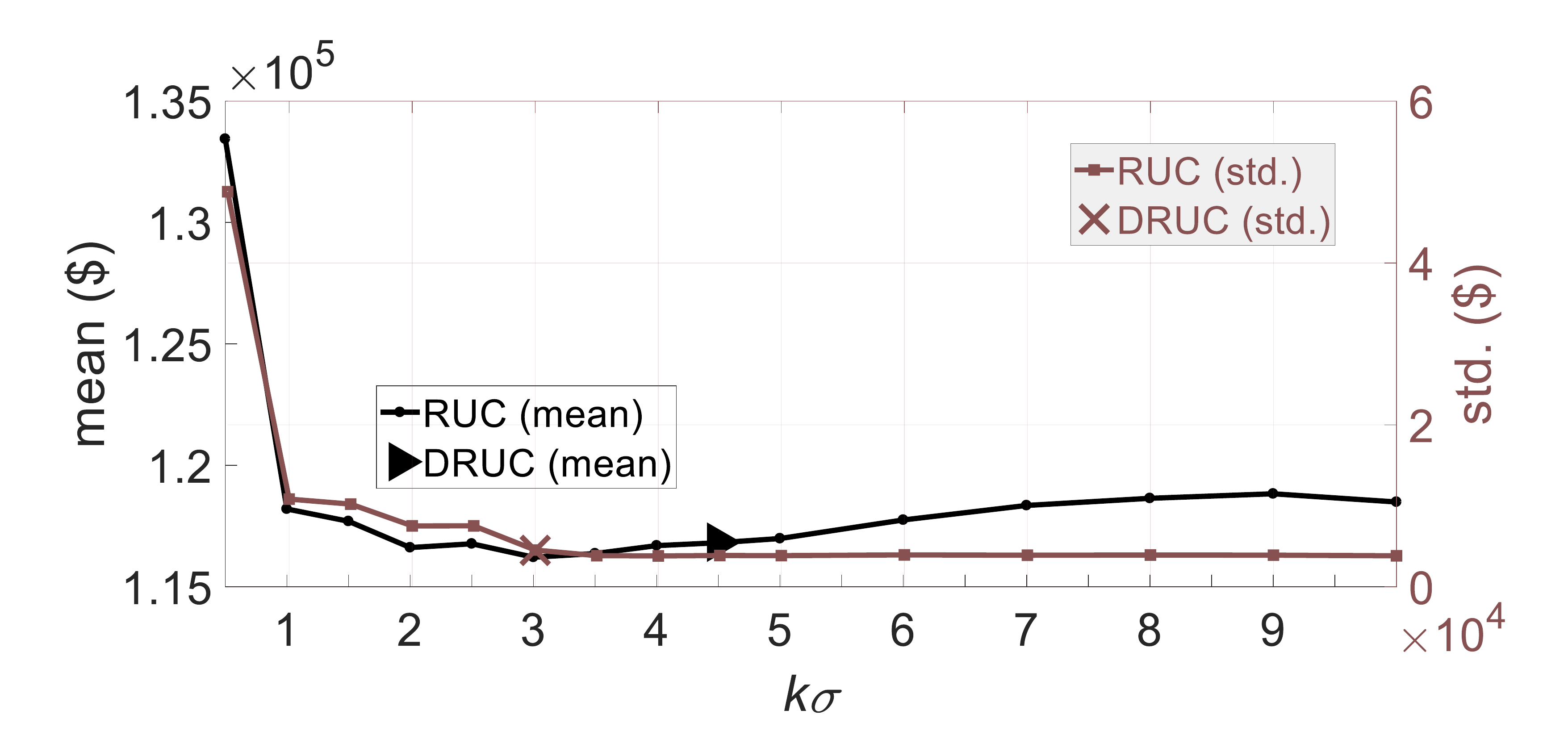}%
\label{fig:multi_k_cost}}
\hfil
\subfloat[]{\includegraphics[width=3.4in]{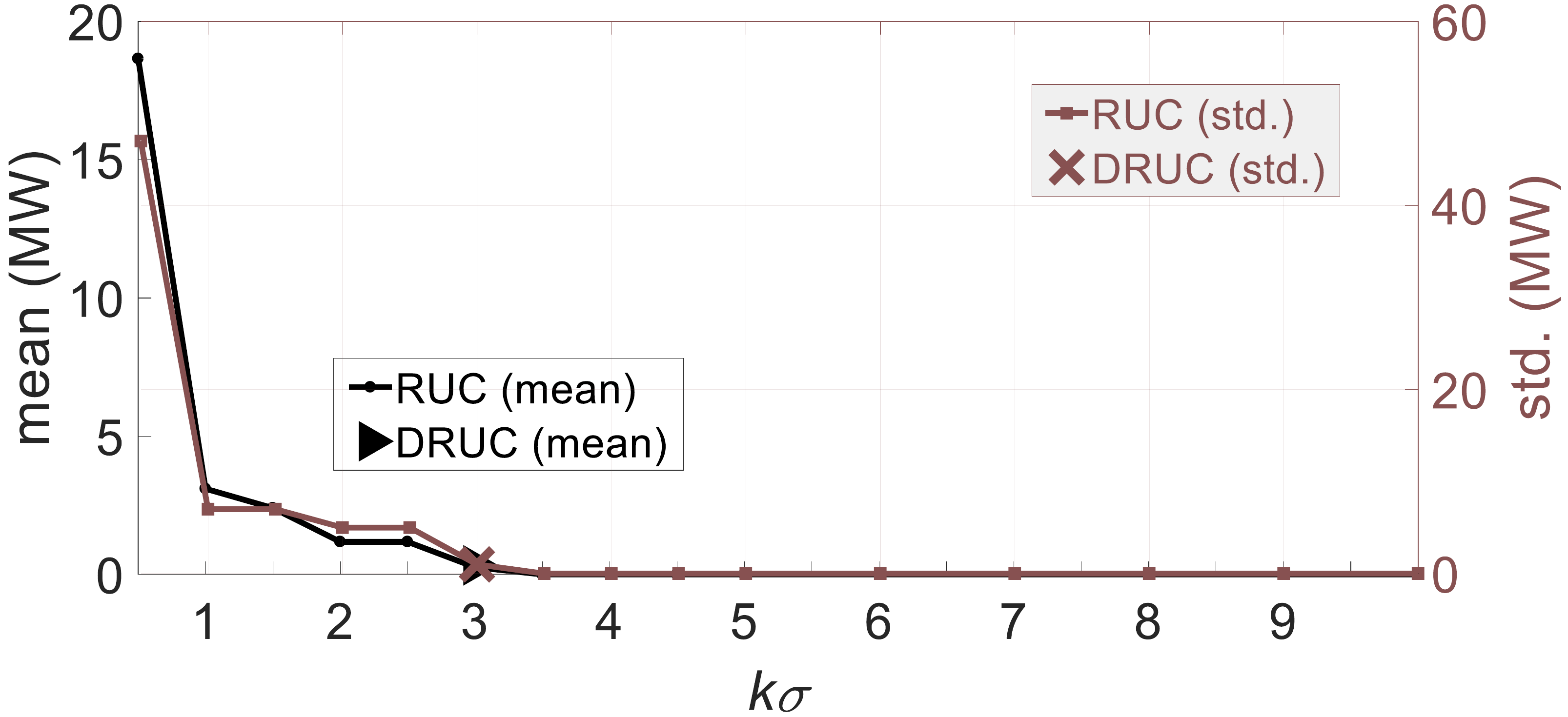}%
\label{fig:multi_k_curtailment}}
\caption{Performances of RUC method with different uncertainty intervals and comparison with the result of DRUC (black axis on the left side for mean; blue axis on the right side for std.). (a) Real cost. (b) Wind curtailment.}
\label{fig:multi_k}
\end{figure}

One may argue that the parameters of RUC model, e.g., the interval width of stochastic parameter and the uncertainty budget, can be tuned to achieve better cost efficiency. This is true as shown in Fig. \ref{fig:multi_k}. In our case, the mean value (or std.) of {\color{black}real-time operation cost} induced by DRUC is equivalent to that of RUC when the uncertainty intervel is about $\pm4.5\bm{\sigma}$ (or $\pm3\bm{\sigma}$). It seems that the $3\bm{\sigma}$-interval is optimal for this case, because the cooresponding {\color{black}real-time operation cost} is minimum, and the wind curtilment amout shown in Fig.~\ref{fig:multi_k}\subref{fig:multi_k_curtailment} is also small enough. However, it is exhuasted to tune the uncertianty set, especially for large-scale systems, since it requires one to pointwisely solve a couple of RUC and ED problems. Another obstacle is that the ¡°optimal parameters¡± are case-dependent, and even sensitive to the operation condition.

\subsection{Computational Efficiency}
\begin{figure}[!t]
\centering
\subfloat[]{\includegraphics[width=3.4in]{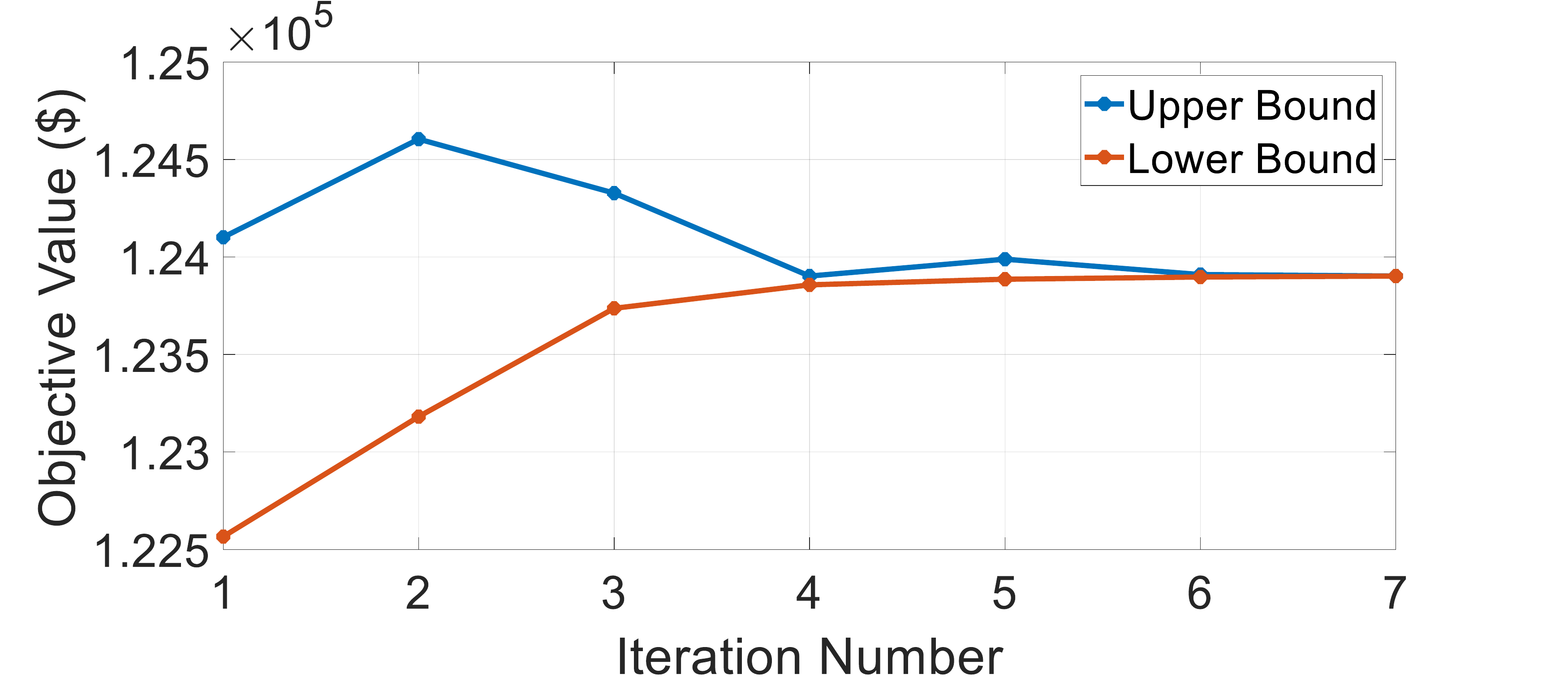}%
\label{fig:convergence_ub_lb}}
\hfil
\subfloat[]{\includegraphics[width=3.4in]{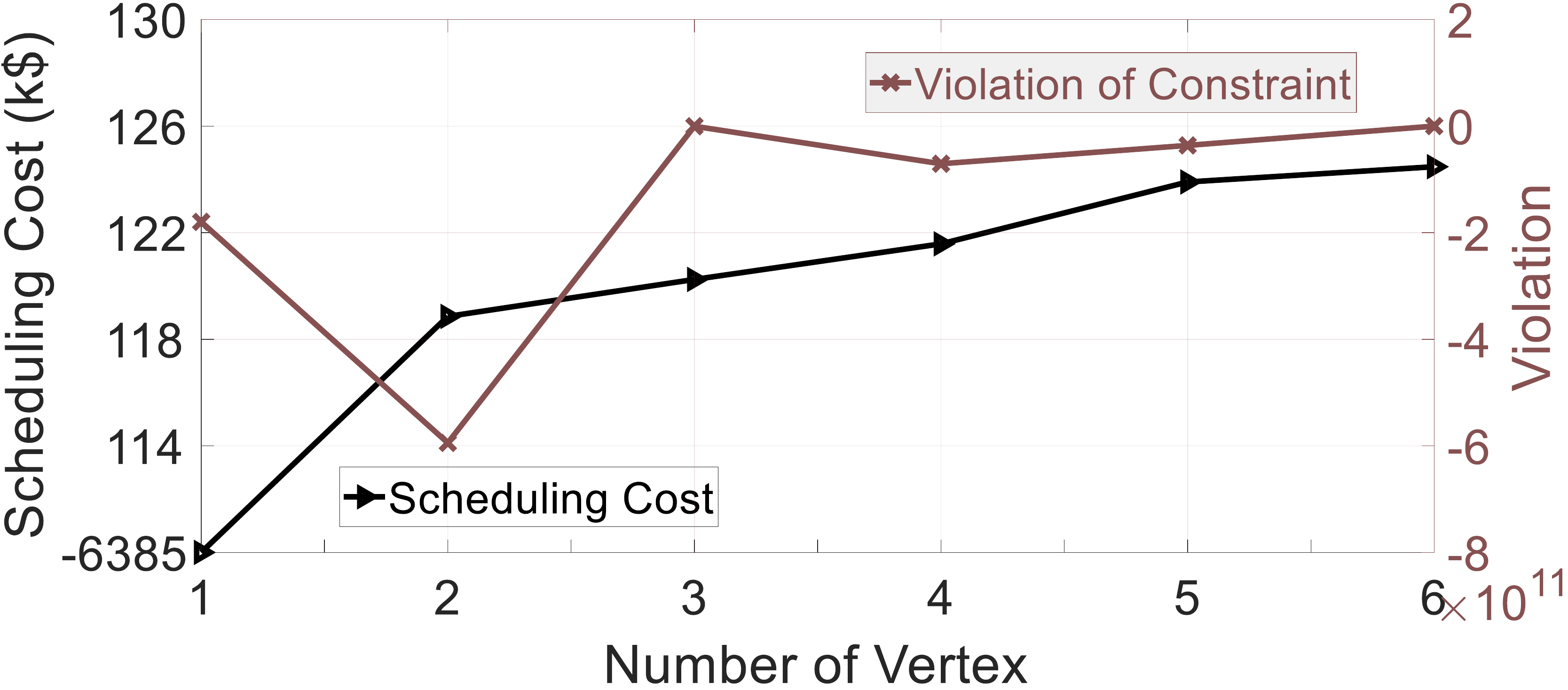}%
\label{fig:convergence_cost_violation}}
\caption{Convergence profiles. (a) Evolution of the upper bound and lower bound of the cutting plane algorithm. (b) Changes of the scheduling cost and the degree of constraint violation as more vertexes are added.}
\label{fig:mean_convergence}
\end{figure}

In this 6-bus system case, it takes 2 to 20 iterations for the cutting plane algorithm to converge, and about 3 to 15 MI-SDP problems should be solved to finish one single DRUC problem. Figure \ref{fig:mean_convergence}\subref{fig:convergence_ub_lb} shows the process of solving one relaxed MI-SDP problem with the cutting plane algorithm, which converges at the 7-th iteration. {\color{black}As more dual cutting planes are added to the master problem, the feasible region gets smaller, and thus the lower bound keeps improving. However, the upper bound derived from the sub-problem may not be monotonously decreasing.}
Figure \ref{fig:mean_convergence}\subref{fig:convergence_cost_violation} shows the monotonically increasing series of scheduling cost as more vertexes are added to the MI-SDP problem. The blue line (right axis) in the same figure demonstrates that the DRUC problem terminates only if the violation degree of constraint (\ref{eqn:Z_dual_cs}) becomes nonnegative, or in other words, the number of vertexes detected and added are sufficient.

Throughout the experiments, we observe that in order to solve the DRUC problem, most of the time are spent on solving the MI-SDP problems. The average runtime of three methods on the 6-bus system case is reported in Table \ref{tab:runtime}.
\begin{table}[t]
\centering
\caption{Comparison of Runtime between UC, RUC and DRUC Methods.\label{tab:runtime}}
\begin{tabular}{c|ccc}
\hline
Method      & UC     & RUC    & DRUC     \\ \hline
Runtime (s) & 0.25 & 0.94 & 271.57 \\ \hline
\end{tabular}
\end{table}

We also employed our algorithm to the DRUC problems of larger systems such as the IEEE 118-bus system, and found that if we omitted the security constraint, some cases can be successfully solved, while most of them run for more than 10 hours and still didn't converge. This is due to the limitations of hardware, the huge number of vertexes in large-scale problems, and more importantly the limited capability of current SDP solvers (the interior point method).

\section{Conclusions}
In our distributionally robust unit commitment model, the second order moment information of stochastic parameters is preserved, and hence the model can capture the correlation of intermittent renewable generations{\color{black}, which cannot be well handled by linear moments. Hence, the decision yielded from our method could be less conservative than that from distributionally robust models with linear moment constraints in general.} The model is parameter free and can still achieve remarkable results, while the robust unit commitment model requires the uncertainty set to be tuned to achieve higher efficiency. Moreover, the set of distributions is adaptive to the data available, i.e., with more data fed into the estimator, the set becomes nearer to the ``true'' one. The only concern about this method is how to improve the computational efficiency, which is also our ongoing research.

We are aware that the MI-SDP reformulation is built on the assumption that the stochastic parameters can take any values in the whole space, which is unnecessary true and may cause the so-call ``tail effect''. However, under the moments constraints, the probability that stochastic parameters take any extreme values should be very low. Still, it is interesting to investigate the exact impact of such an assumption. From another perspective, since we are able to derive the empirical distribution (by using other statistics method like Bayesian nonparametric models, even the exact formulation of the distribution can be inferred), we can then employ this distribution to a traditional stochastic programming. The comparison of the abovementioned method and distributionally robust optimization is also an attractive research direction.



%

\appendices
\section{Unit Commitment Model}\label{append:UC}
In the (distributionally) robust models, reserve requirements is neglected. Besides, nodal injections of RES are integrated, resulting in the unit commitment model as follows,
\begin{subequations} \label{eqn:SCUC_full}
    \begin{align}
        \nonumber & \underset{\bm{x},\bm{y},\bm{\xi}}{\mathop{\min }}\, \sum\limits_{g\in {{\mathsf{\mathcal{G}}}}}{\sum\limits_{t\in \mathsf{\mathcal{T}}}{{{x}_{g,t}}N{{L}_{g,t}}+{{u}_{g,t}}S{{U}_{g,t}}}} \\
        & ~~~~~~~~~~~~~~~~~+{{v}_{g,t}}S{{D}_{g,t}}+{{C}_{g,t}}\left( {{p}_{g,t}} \right) \\
        \nonumber & \operatorname{s.t.} \\
        & {{x}_{g,t}}-{{x}_{g,t-1}}={{u}_{g,t}}-{{v}_{g,t}} ~~\forall g\in\mathcal{G},t\in\mathcal{T} \label{eqn:SCUC_full_status}\\
        & \sum\limits_{\tau =\max \left\{ 1,t-M{{U}_{g}}+1 \right\}}^{t}{{{u}_{g,\tau }}}\le {{x}_{g,t}} ~~\forall g\in\mathcal{G},t\in\mathcal{T} \label{eqn:SCUC_full_min_on}\\
        & \sum\limits_{\tau =\max \left\{ 1,t-M{{D}_{g}}+1 \right\}}^{t}{{{v}_{g,\tau }}}\le 1-{{x}_{g,t}} ~~\forall g\in\mathcal{G},t\in\mathcal{T} \label{eqn:SCUC_full_min_down}\\
        & P^{\min}_{g}{{x}_{g}} \le {{p}_{g,t}} \le P^{\max}_{g}{{x}_{g}} ~~\forall g\in\mathcal{G},t\in\mathcal{T} \label{eqn:SCUC_full_gen_limit} \\
        & -R_{g}^{-}\le {{p}_{g,t}}-{{p}_{g,t-1}}\le R_{g}^{{+}} ~~\forall g\in {{\mathsf{\mathcal{G}}}},t,t-1\in \mathsf{\mathcal{T}} \label{eqn:SCUC_full_ramping_limit}\\
        & \theta_{n,t} = 0 ~~ n\in\mathcal{N}_{\rm{ref}} \label{eqn:SCUC_full_ref} \\
        & |(\theta_{n,t}-\theta_{m,t})/X_{(m,n)}| \le F_{(m,n)} ~~\forall (m,n)\in\mathcal{L} \label{eqn:SCUC_full_flow_limit}\\
        \nonumber & \sum_{m\in\mathcal{N}}B_{n,m}\theta_{m,t} = \sum_{g\in\mathcal{G}}\kappa_{g,n}^{\mathcal{GN}}p_{g,t} - D_{n,t} \\
        & ~~~~~~~~~~~~~~~~~~~~~+ \sum_{e\in\mathcal{E}}\kappa_{e,n}^{\mathcal{EN}}\xi_{e,t} ~~\forall n\in\mathcal{N},t\in\mathcal{T} \label{eqn:SCUC_full_dc_equation}
    \end{align}
\end{subequations}
\begin{itemize}
  \item $NL_{g,t},~SU_{g,t},~SD_{g,t},~C_{g,t}(\cdot)$:~No-load, start-up, shut-down costs and variable cost function of generator $g$ at time $t$.
  \item $MU_{g},~MD_{g}$:~Minimum-up and minimum-down time of unit $g$.
  \item $R^{+}_{g},~R^{-}_{g}$:~Ramp-down and ramp-up limits of unit $g$.
  \item $X_{(m,n)},~F_{(m,n)}$:~Reactance and rating of transmission line connecting buses $m$ and $n$.
  \item $B_{m,n}$:~Element on the intersection of the $m$-row and the $n$-th column of the nodal susceptance matrix.
  \item $\kappa_{g,n}^{\mathcal{GN}},~\kappa_{e,n}^{\mathcal{EN}}$:~0-1 coefficient indicating whether unit $g$ or RES node $e$ located at bus $n$.
  \item $D_{n,t}$:~Load demand at bus $n$ at time $t$.
  \item $P^{\min}_{g},~P^{\max}_{g}$:~Minimum and maximum production levels of unit $g$.
  \item $\mathcal{T},~\mathcal{G},~\mathcal{L},~\mathcal{N},~\mathcal{E},~\mathcal{N}_{\rm{ref}}$:~Sets of time periods, units, transmission lines (represented by the two ends), buses, RES nodes and node of reference bus.
  \item $x_{g,t},~u_{g,t},~v_{g,t}~(\bm{x})$:~Binary variables indicating whether the unit is on, started up and shut down.
  \item $p_{g,t},~\theta_{n,t}~(\bm{y})$:~Production level of unit $g$ and phase angel of bus $n$ at time $t$.
  \item $\xi_{e,t}~(\bm{\xi})$:~Output of RES node $e$ at time $t$.
\end{itemize}

Constraints~(\ref{eqn:SCUC_full_status})-(\ref{eqn:SCUC_full_min_down}) include state transition equations of units and minimum up/down time limits of units. Constraints~(\ref{eqn:SCUC_full_gen_limit})-(\ref{eqn:SCUC_full_dc_equation}) are production and ramping limits of units, power flow limits on transmission lines, the condition of reference bus and DC power flow equation respectively.

{\color{black}\section{Dualization of Semi-infinite Programming}\label{append:dualization}
First, let us define the \emph{Lagrangian} associated the problem~(\ref{eqn:Z_primal}) as,
\begin{align*}
& \mathcal{L}(f_{\bm{\xi}}\ge0,h_0,\bm{h},\bm{H}) = \int_{\Xi}{Q(\bm{x,\xi} ){{f}_{\bm{\xi}}}\mathrm{d}\bm{\xi} }  \\
& \qquad - h_0\left(\int_{\Xi}{{{f}_{\bm{\xi}}}\mathrm{d}\bm{\xi}-1}\right) - \sum_{i=1}^{|\bm{\xi}|} {h}_{i}\left(\int_{\Xi}{{{\xi}_{i}}{{f}_{\bm{\xi}}}\mathrm{d}\bm{\xi}-{{\bar{\xi}}_{i}}}\right) \\
& \qquad \qquad  - \sum_{j=1,k=1}^{|\bm{\xi}|,|\bm{\xi}|}H_{j,k}\left(\int_{\Xi }{{{\xi}_{j}}{{\xi }_{k}}{{f}_{\bm{\xi}}}\mathrm{d}\bm{\xi}-{{\Sigma}_{jk}}+{{\bar{ \xi }}_{j}}{{\bar{\xi}}_{k}}}\right).
\end{align*}
The \emph{Lagrange dual function} that provides an upper bound for the problem~(\ref{eqn:Z_primal}), is given as,
\begin{align*}
g(h_0,\bm{h},\bm{H}) = \sup_{f_{\bm{\xi}}\ge0}\, \mathcal{L}(f_{\bm{\xi}},h_0,\bm{h},\bm{H}).
\end{align*}
The parameters associated with $f_{\bm{\xi}}$ should be non-positive, since otherwise $g(h_0,\bm{h},\bm{H})$ would be $+\infty$.
Therefore, the \emph{Lagrange dual problem}, which aims to minimizes the upper bound, can be formulated as follows,
\begin{align*}
& \min_{h_0,\bm{h},\bm{H}}\, g(h_0,\bm{h},\bm{H})\\
& \overset{g<+\infty}\Longleftrightarrow  \\
& \min_{h_0,\bm{h},\bm{H}}\, h_0+\bm{h}^{\top}\bar{\bm{\xi}}+\bm{H}(\bm{\Sigma}+\bar{\bm{\xi}}\bar{\bm{\xi}}^{\top})\\
& {~~\rm{s.t.~}}  Q(\bm{x,\xi}) - ({{h}_{0}} + {{\bm{h}}^{\top}}\bm{\xi} + {{\bm{\xi} }^{\top}}\bm{H\xi}) \le 0 ~~ \forall \bm{\xi} \in \Xi.
\end{align*}
The problem above can be further written as the form of problem~(\ref{eqn:Z_dual}).


}
\section*{Acknowledgment}
The first author would like to thank Johan L{\"o}fberg for the discussion and suggestion on solving the MI-SDP problem, who also develops and maintains the YALMIP toolbox for modeling and optimization in MATLAB.

\ifCLASSOPTIONcaptionsoff
  \newpage
\fi



\bibliographystyle{IEEEtran}
\bibliography{IEEEabrv,DR_MISDP}
\end{document}